\documentclass[11pt]{article}
\usepackage{latexsym}
\usepackage{amssymb}
\newlength{\mytopmargin}
\newlength{\myleftmargin}
\newtheorem{lemma}{Lemma}[section]

\newtheorem{thm}[lemma]{Theorem}
\newtheorem{cor}[lemma]{Corollary}
\newtheorem{prop}[lemma]{Proposition}
\newcommand{\nd}{\mathcal{D}}
\newcommand{\cc}{\mathbb{C}}

\newcommand{\bin}[2]{({ {\scriptscriptstyle #1} \atop
{\scriptscriptstyle #2} })}
\newcommand{\inn}[2]{\left\langle #1, #2 \right\rangle}

\begin{document}
\noindent
\begin{center}{ \Large\bf Isomorphisms of type $A$ affine Hecke algebras
\\[1.5mm] and multivariable orthogonal polynomials}
\end{center}
\vspace{5mm}
\noindent
\begin{center} T.~H.~Baker 
and P.~J.~Forrester \\[3mm]
{\it Department of Mathematics, University of Melbourne, \\
Parkville, Victoria 3052, Australia} 
\end{center}
\vspace{.5cm}

\begin{quote} We examine two isomorphisms between affine Hecke algebras
of type $A$ associated with parameters $q^{-1}$, $t^{-1}$  and
$q$, $t$. One of them maps the non-symmetric Macdonald polynomials
$E_{\eta}(x;q^{-1},t^{-1})$ onto $E_{\eta}(x;q,t)$, while the other
maps them onto non-symmetric analogues of the multivariable Al-Salam \& Carlitz
polynomials. Using the properties of
$E_{\eta}(x;q^{-1},t^{-1})$, the corresponding properties of these
latter polynomials can then be elucidated.

\end{quote}
 
\section{Introduction}

In several recent works \cite{uji96a}--\cite{uji97a},
\cite{kakei96}--\cite{kakei97a},
eigenstates of the rational (type $A$) Calogero-Sutherland model
have been investigated from an algebraic point of view. In particular
it has been shown that the algebra
governing the eigenfunctions of the {\it periodic} Calogero-Sutherland
model (namely the type $A$ degenerate affine Hecke algebra augmented
by type $A$ Dunkl operators) is isomorphic to its {\it rational} model
counterpart. This enables information to be gleaned about the
properties of the eigenfunctions in the rational case (the (non-)symmetric
Hermite polynomials) from the
corresponding periodic eigenfunctions (the (non-)symmetric
Jack polynomials). 

To summarize the argument, consider the type $A$ Dunkl operators
$$
d_i := \frac{\partial}{\partial x_i}
+ \frac{1}{\alpha} \sum_{p\neq i} \frac{1-s_{ip}}{x_i-x_p}
$$
which, along with the operators representing multiplication by the
variable $x_i$ and the elementary transpositions $s_{ij}$, satisfy 
the following commutation relations
\begin{eqnarray}
{}[d_i, x_j] &=& \left\{ \begin{array}{cc}
-\frac{1}{\alpha} s_{ij} & i\neq j\\[2mm]
1 + \frac{1}{\alpha} \sum_{p\neq i} s_{ip} & i=j
\end{array} \right. \nonumber\\
d_i\,s_{ip} &=& s_{ip}\,d_p \hspace{3cm}
[d_i, s_{jp}] =0, \qquad i\neq j,p \label{dunklalg}
\end{eqnarray}
It is easily checked that the map $\rho$ defined by
\begin{equation} \label{iso.kakei}
\rho(x_i) = \frac{1}{2}(x_i - d_i), \qquad \rho(d_i) = d_i, \qquad
\rho(s_{ij}) = s_{ij}
\end{equation}
is an isomorphism of the algebra (\ref{dunklalg}) \cite{uji96a}.

Now, the non-symmetric Jack polynomials $E_{\eta}(x)$, indexed by
compositions $\eta:=(\eta_1,\ldots,\eta_n)$ can be defined \cite{opdam95a} 
as the
unique eigenfunctions of the mutually commuting Cherednik operators
\begin{equation} \label{cher.1}
\xi_i := \alpha x_i d_i + \sum_{p>i} s_{ip} - n + 1
\end{equation}
with a unique expansion of the form
\begin{equation} \label{expansion}
E_{\eta}(x) = x^{\eta} + \sum_{\nu<\eta} c_{\eta\nu}\, x^{\nu} .
\end{equation}
Here, the partial order $<$ is defined on compositions by:
$\nu < \eta$ iff $\nu^+ < \eta^+$ with respect to
the dominance order (where $\nu^+$ is the unique
partition associated to $\nu$ etc) {\it or} $\nu^+ = \eta^+$,
$\nu \ne \eta$ and
$\sum_{i=1}^p (\eta_i - \nu_i) \geq 0$, for all $p=1,\ldots n$.
The polynomial $E_{\eta}(x)$ is an eigenfunction of $\xi_i$ given
by (\ref{cher.1}) with eigenvalue
\begin{equation} \label{e-val}
\bar{\eta}_i = \alpha\eta_i - \#\{k<i | \eta_k \geq \eta_i \} -
\#\{k>i | \eta_k > \eta_i \}
\end{equation}

Using the isomorphism (\ref{iso.kakei}) it follows that the polynomials
\cite{takemura97a,roesler97a,kakei97a}
$$
E^{(H)}_{\eta}(x) := E_{\eta}(\rho(x))\,.\,1
$$ 
are eigenfunctions of the operators
\begin{equation} \label{hermite.op}
h_i = \rho(\xi_i) = \xi_i - \frac{\alpha}{2} d_i^2
\end{equation}
which are precisely the eigenoperators of the non-symmetric
Hermite polynomials \cite{forr96c}. The orthogonality of these latter
polynomials with respect to the usual multivariable Hermite
inner product then follows from the fact that
the operator (\ref{hermite.op}) is self-adjoint with respect to the
inner product
\begin{equation}
\inn{f}{g} := \prod_{i=1}^n \int_{-\infty}^\infty
dx_i \, e^{-x_i^2} \prod_{1 \le j < k \le n} |x_j - x_k|^{2/\alpha}
\,f\, g
\label{inn.h}
\end{equation}

In this work, we provide a similar analysis of the Macdonald case.
As such, we introduce an isomorphism of the $q$-analogue of the
algebra (\ref{dunklalg}), namely the (type $A$) affine Hecke algebra
augmented by $q$-Dunkl operators. To describe this mapping, we need
to introduce some further concepts. 

The generalization of the formalism of non-symmetric Jack polynomials to
the Macdonald case involves replacing the Cherednik operators (\ref{cher.1})
by their $q$-analogues which can be constructed by means of the 
generators of the affine Hecke algebra \cite{mac95}. In the type $A$
case, one can describe this using the Demazure-Lustig operators
\begin{eqnarray} \label{maldonado}
T_i &:=& t + \frac{tx_i-x_{i+1}}{x_i - x_{i+1}}\left(s_i -1 \right)
\hspace{2cm} i=1,\ldots,n-1 \\
T_0 &=& t + \frac{qtx_n-x_1}{qx_n - x_1}\left(s_0 -1 \right)
\end{eqnarray}
along with the operator
\begin{equation}\label{omega}
\omega := s_{n-1}\cdots s_2\,s_1\tau_1 =
s_{n-1}\cdots s_i\tau_i s_{i-1}\cdots s_1, 
\end{equation}
where $\tau_i$ is the operator which replaces $x_i$ by $qx_i$.
The affine Hecke algebra is then generated by elements $T_i$,
$0\leq i\leq n-1$ and $\omega$, satisfying the relations
\begin{eqnarray}
(T_i-t)\,(T_i+1) &=& 0 \label{tdefs.1}\\
T_i\;T_{i+1}\;T_i &=& T_{i+1}\;T_i\;T_{i+1} \label{tdefs.2}\\
T_i\;T_j &=& T_j\;T_i \qquad |i-j| \geq 2 \\
\omega\;T_i &=& T_{i-1}\;\omega  \label{tdefs.4}
\end{eqnarray}
There is a commutative subalgebra generated by elements of the form
\cite{cher91a,cher94a}
\begin{equation}\label{yi}
Y_i := t^{-n+i}\, T_i\cdots T_{n-1}\,\omega\,T_1^{-1}\cdots T_{i-1}^{-1}
\end{equation}
which have the following relations with the generators $T_i$
\begin{equation}\label{ty}
T_i\,Y_{i+1} = t Y_i\,T_i^{-1}, \qquad
T_i\,Y_i = Y_{i+1}\,T_i + (t-1)Y_i, \qquad [T_i,Y_j]=0, \: \:
j \ne i,i+1.
\end{equation}
The non-symmetric Macdonald polynomials $E_{\eta}(x;q,t)$ are defined
as the simultaneous eigenfunctions of the commuting operators $Y_i$ with
an expansion of the form (\ref{expansion}). The corresponding eigenvalue
is $t^{\bar{\eta}_i}$ with $\bar{\eta}_i$ given in (\ref{e-val}).
{}From now on, we drop the dependence on $q$ and $t$ and just write
$E_{\eta}(x)\equiv E_{\eta}(x;q,t)$ when the meaning is unambiguous.

Define the following degree-raising operator
\begin{equation} \label{def.e}
e_i = t^{i-1} T_i\cdots T_{n-1}\,x_n\,\omega\,T_1^{-1} \cdots
T_{i-1}^{-1} .
\end{equation}
Using (\ref{tdefs.2})--(\ref{tdefs.4})
it can be shown that the operators $e_i$ form a set of mutually
commuting  operators. Our first result is
\begin{thm} \label{nata}
We have
$$
E_{\eta}(e_1,\ldots,e_n;q^{-1},t^{-1}) \,.\,1 = \alpha_{\eta}(q,t)\,
E_{\eta}(x_1,\ldots,x_n;q,t)
$$
where
\begin{equation} \label{cuerda}
\alpha_{\eta}(q,t) = q^{\sum_i \bin{\eta_i}{2}} t^{\sum_i(n-i)\eta_i^+
-\ell(w_{\eta})}
\end{equation}
with $\ell(w_{\eta})$ the length of the (unique) minimal permutation
sending $\eta$ to $\eta^+$.
\end{thm}

The symmetric Al-Salam \& Carlitz (ASC) polynomials were examined in
\cite{forr97c} as $q$-analogues of multivariable Hermite polynomials.
There are two families of ASC polynomials, denoted $U_\lambda^{(a)}(x;q,t)$
and $V_\lambda^{(a)}(x;q,t)$, which are simply related by
\begin{equation}\label{UV}
V_\lambda^{(a)}(x;q^{-1},t^{-1}) = U_\lambda^{(a)}(x;q,t).
\end{equation}
The polynomials $V_\lambda^{(a)}$
can be defined as the unique polynomials of the form 
$$
V^{(a)}_{\lambda}(x;q,t) = P_{\lambda}(x;q,t) + \sum_{\mu<\lambda} 
b_{\lambda\mu} P_{\mu}(x;q,t)
$$
which are orthogonal with respect to the inner product
\begin{equation}
\label{innerv}
\langle f , g \rangle^{(V)} :=
\int_{[1,\infty]^n} f(x)g(x) \, d_q\mu^{(V)}(x), \quad
d_q\mu^{(V)}(x) :=  \Delta_q^{(k)}(x)
\prod_{l=1}^n w_V(x_l;q) d_qx_l .
\end{equation}
Here, $P_{\lambda}(x;q,t)$ denotes the symmetric Macdonald polynomial 
\cite{mac} and
we use the notation for $q$-integrals
\begin{equation}
\int_1^\infty f(x) d_qx := (1-q) \sum_{n=0}^\infty f(q^{-n})
q^{-n}
\end{equation}
while
\begin{eqnarray}
w_V(x;q) &=& {(q ;q)_\infty ({1 \over a};q)_\infty
(qa;q)_\infty \over {(x;q)'}_\infty
({x \over a};q)_\infty} \nonumber \\
\Delta_q^{(k)}(x_1,\ldots,x_n) &:=& \prod_{p=-(k-1)}^k
\prod_{1 \le i < j \le n}(x_i - q^px_j) \label{delta},
\end{eqnarray}
where the dash in $(x;q)_\infty'$ denotes that any vanishing factor 
is to be deleted, and it is assumed $a<0$.

The polynomials $U_\lambda^{(a)}$ are orthogonal with respect to the
inner product
\begin{equation}\label{inneru}
\langle f | g \rangle^{(U)} :=
\int_{[a,1]^n} f(x)g(x) \, d_q\mu^{(U)}(x), \quad
d_q\mu^{(U)}(x) :=   \Delta_q^{(k)}(x)
\prod_{l=1}^n w_U(x_l;q) d_qx_l
\end{equation}
where $\Delta_q^{(k)}$ is given by (\ref{delta}) and
\begin{equation}\label{wu}
 w_U^{(a)}(x;q) := {(qx;q)_\infty ({qx \over a};q)_\infty \over
 (q;q)_\infty (a;q)_\infty ({q \over a};q)_\infty}
 \end{equation}
\begin{equation}\label{defqi}
\int_a^1 f(x) \, d_qx := (1-q) \left(
\sum_{n=0}^\infty f(q^n) q^n - a\sum_{n=0}^\infty f(aq^n) q^n \right),
\quad (a <0)
\end{equation}
This can be regarded as a consequence of (\ref{UV}), and the formulas
\begin{equation}\label{qiuv}
{1 \over 1 - q}
\int_a^1 w_U^{(a)}(x;q) f(x) \, d_qx \Big |_{q \mapsto q^{-1}} =
{1 \over 1 - q} \int_1^\infty w_V^{(a)}(x;q) f(x) \, d_qx
\end{equation}
\begin{equation}\label{deltain}
\Delta_{q^{-1}}^{(k)}(x) = q^{-kn(n-1)} \Delta_{q}^{(k)}(x^R)
\end{equation}
where $x^R = (x_n,x_{n-1},\dots,x_1)$. The formula (\ref{qiuv}) is
established in \cite[eq.~(2.23)]{forr97c}, while (\ref{deltain})
follows immediately from the definition (\ref{delta}).

Non-symmetric analogues of the ASC polynomials can be introduced in
the following manner: consider the following $q$-analogues of the 
type $A$ Dunkl operators \cite{dunkl89a} examined in \cite{forr97b},
\begin{equation} \label{dunk.old}
D_i = x_i^{-1} \left( 1-t^{n-1}\,T_i^{-1}\cdots T_{n-1}^{-1}\,\omega\,
T_1^{-1}\cdots T_{i-1}^{-1} \right)
\end{equation}
and let 
\begin{equation}
E_i := D_i + (1+a^{-1})t^{n-1}\,Y_i - a^{-1}\,e_i \label{e.def}.
\end{equation}
The operators $E_i$ mutually commute, and 
our second main result is that
\begin{thm} \label{yanan}
The polynomials
\begin{equation} \label{non.v}
E^{(V)}_{\eta}(x;q,t) = \frac{(-a)^{|\eta|}}{\alpha_{\eta}(q,t)}\:
E_{\eta}(E;q^{-1},t^{-1})\,.\,1
\end{equation}
where $\alpha_{\eta}(q,t)$ is given by (\ref{cuerda})
are the unique polynomials with an expansion of the form
$$
E^{(V)}_{\eta}(x;q,t) =  E_{\eta}(x;q,t) + \sum_{|\nu| < |\eta|}
c_{\eta\nu}\, E_{\nu}(x;q,t)
$$
which are orthogonal with respect to the inner product (\ref{innerv}).
Furthermore, these polynomials are simultaneous eigenfunctions of
the commuting family of eigenoperators
\begin{equation}\label{def.hop}
h_i = Y_i + (1+a) t^{1-n} D_i + a t^{2-2n} D_i Y_i^{-1} D_i
\end{equation}
with eigenvalue $t^{\bar{\eta}_i}$.
\end{thm}

An immediate consequence of Thm.~\ref{yanan}, (\ref{UV}), and (\ref{qiuv}),
(\ref{deltain}) is

\begin{cor} \label{cor1.3} The polynomials
\begin{equation}\label{uv}
E_\eta^{(U)}(x;q,t) := E_\eta^{(V)}(x^R;q^{-1},t^{-1})
\end{equation}
are the unique polynomials with an expansion of the form
$$
E^{(U)}_{\eta}(x;q,t) =  E_{\eta}(x^R;q^{-1},t^{-1}
) + \sum_{|\nu| < |\eta|}
d_{\eta\nu}\, E_{\nu}(x^R;q^{-1},t^{-1})
$$
which are orthogonal with respect to the inner product (\ref{inneru}).
These polynomials are simultaneous eigenfunctions of the operators
$\hat{h}_i$, where $\hat{h}_i$ denotes the operator (\ref{def.hop}) modified by
the involution $\hat{\,\,}$, which is defined by the mappings $q \mapsto q^{-1}$,
$t \mapsto t^{-1}$ and $x \mapsto x_{n+1-i}$.
\end{cor}

In section 2, we examine the various properties of non-symmetric Macdonald
polynomials used in subsequent calculations, including raising and
lowering operators, and introduce a non-symmetric analogue of Kaneko's
kernel \cite{kaneko93a}. We finish the section with a proof of
Thm.~\ref{nata}. An isomorphism between Hecke algebras is introduced in
section 3, facilitating a proof of Thm.~\ref{yanan}. Various
properties of these non-symmetric ASC polynomials are then described
including their normalization and a generating function. We conclude
by clarifying their relationship to the non-symmetric analogues of the
shifted Macdonald polynomials.

\setcounter{equation}{0}
\section{Non-symmetric Macdonald polynomials}

In this section we gather together some (old and new) results concerning
non-symmetric Macdonald polynomials $E_{\eta}(x)$ in preparation of
the proof of Thm.~\ref{nata}, as well as the forthcoming section on the
non-symmetric ASC polynomials.

For future reference we note that the operators $T_i$ and $\omega$ 
defined by (\ref{maldonado}) and (\ref{omega}) have the properties
\begin{eqnarray}\label{2.1}
T_i^{-1}\,x_{i+1} =t^{-1}x_i\,T_i \hspace*{2cm} &
T_i^{-1}\,x_i = x_{i+1}\,T_i^{-1} + (t^{-1}-1)x_i  \nonumber\\
T_i\,x_i = t x_{i+1}\,T_i^{-1} \hspace*{2cm} &
T_i\,x_{i+1} = x_i\,T_i +(t-1)x_{i+1} \label{id.2}\\
\omega\,x_1 = qx_n\omega \hspace{2cm} & \omega\,x_{i+1} = x_i\,\omega
\nonumber
\end{eqnarray}
valid for $1\leq i\leq n-1$.
Also note the following action of $T_i$ on monomials
\begin{equation} \label{t.action}
T_i\, x_i^a x_{i+1}^b = \left\{ \begin{array}{ll}
(1-t)x_i^{a-1}x_{i+1}^{b+1} +
\cdots +(1-t)x_i^{b+1}x_{i+1}^{a-1} + x_i^bx_{i+1}^a& a > b \\
t x_i^ax_{i+1}^a & a=b \\
(t-1)x_i^{a}x_{i+1}^{b} +
\cdots +(t-1)x_i^{b-1}x_{i+1}^{a+1} + t x_i^bx_{i+1}^a & a < b
\end{array} \right.
\end{equation}

There exists a variant of the $q$-Dunkl operator (\ref{dunk.old}) 
which is relevant to the forthcoming discussion. With $\hat{\,\,}$ denoting
the involution defined in the statement of Corollary \ref{cor1.3}, this
operator is defined as
\begin{eqnarray}\label{dunk.new}
\nd_i& := & -q \hat{D}_{n+1-i}  \nonumber \\
& = & - qx_i^{-1} \Big ( 1 - t^{-n+1} T_{i-1} \cdots T_1 \omega^{-1}
T_{n-1} \cdots T_i \Big ) \nonumber \\
& = & q t^{-2n + i + 1} D_i Y_i^{-1} T_i \cdots T_{n-1} T_{n-1} \cdots T_i
\end{eqnarray}
In obtaining the first equality in (\ref{dunk.new}), the facts that
\begin{equation}\label{hattw}
\hat{T}_i = T_{n-i}^{-1} \quad \mbox{and} \quad \hat{\omega} = \omega^{-1}
\end{equation}
have been used in applying the operation $\hat{\,\,}$ to (\ref{dunk.old}),
while the second equality can be verified by substituting for
$ Y_i^{-1}$ using (\ref{yi}) and for $D_i$ using (\ref{dunk.old})
and comparing with the first equality.

Since the $D_i$ commute, it follows from the definition of $\nd_i$ that
the $\{ \nd_i \}$ also form a commuting set. Moreover, using (\ref{hattw}),
one can check that the operators $\nd_i$ possess the same 
relations with the generators $T_i$, $\omega$ as do the $D_i$, namely
\begin{eqnarray}
T_i \nd_{i+1} &=& t \nd_i T_i^{-1}, \qquad T_i\nd_i = \nd_{i+1} T_i +
(t-1)\nd_i, \quad 1\leq i\leq n-1 \nonumber\\
{}[T_i, \nd_j] &=& 0, \hspace{7cm} j\neq i, i+1 \label{nd.id}\\
\nd_n \,\omega &=& q\,\omega\, \nd_1 , \hspace{2cm}
\nd_i\, \omega = \omega \,\nd_{i+1} \hspace{2cm} 1\leq i\leq n-1 \nonumber
\end{eqnarray}

To conclude the preliminaries, we follow
Sahi \cite{sahi96a} and introduce the generalized arm and
leg (co-)lengths for a node $s\in\eta$ via
\begin{eqnarray}
a(s)= \eta_i - j && l(s) = \#\{k>i|j\leq \eta_k\leq\eta_i\} \;+\;
\#\{k<i|j\leq \eta_k+1\leq\eta_i\} \nonumber\\
a'(s)=j - 1 && l'(s) = \#\{k>i| \eta_k > \eta_i\} \;+\;
\#\{k<i|\eta_k\geq\eta_i\}  \label{guion}
\end{eqnarray}
and define the quantities
\begin{eqnarray}\label{qu}
d_{\eta}(q,t) &:=& \prod_{s\in\eta} \left( 1-q^{a(s)+1}t^{l(s)+1} \right)
\hspace{2cm} l(\eta) := \sum_{s\in\eta} l(s) \nonumber\\
d'_{\eta}(q,t) &:=& \prod_{s\in\eta} \left( 1-q^{a(s)+1}t^{l(s)} \right)
\hspace{2cm} l'(\eta):= \sum_{s\in\eta} l'(s) \label{constants}\\
e_{\eta}(q,t) &:=& \prod_{s\in\eta} \left( 1-q^{a'(s)+1}t^{n-l'(s)} \right)
\hspace{2cm} a(\eta) := \sum_{s\in\eta} a(s)  \nonumber
\end{eqnarray}
The statistics $l(\eta)$, $l'(\eta)$ and $a(\eta)$ generalize the 
quantity 
\begin{equation}\label{def.n}
b(\lambda):=\sum_i (i-1)\lambda_i = \sum_i \left ( {\lambda'_i\atop 2}\right )
\end{equation}
from partitions to compositions.
{}From \cite{sahi96a} these quantities have the following properties
\begin{lemma} \label{props}
Let $\Phi\eta:=(\eta_2,\ldots,\eta_n,\eta_1+1)$. We have
\begin{eqnarray*}
\frac{d_{\Phi\eta}(q,t)}{d_{\eta}(q,t)} &=&
\frac{e_{\Phi\eta}(q,t)}{e_{\eta}(q,t)} = 1-qt^{n+\bar{\eta}_1},
\quad
\frac{d'_{\Phi\eta}(q,t)}{d'_{\eta}(q,t)} = 1-qt^{n-1+\bar{\eta}_1}, \quad
e_{s_i\eta}(q,t) = e_{\eta}(q,t), \\
\frac{d_{s_i\eta}(q,t)}{d_{\eta}(q,t)} &=& \frac{1-t^{\delta_{i,\eta}+1}}
{1-t^{\delta_{i,\eta}}}, \qquad 
\frac{d'_{s_i\eta}(q,t)}{d'_{\eta}(q,t)} = \frac{1-t^{\delta_{i,\eta}}}
{1-t^{\delta_{i,\eta}-1}} \qquad \mbox{{\rm for} $\eta_i > \eta_{i+1}$} \\
a(\Phi\eta) &=& \eta_1 + a(\eta),  \qquad l(\Phi\eta) = l(\eta)
+ \#\{k>1| \eta_k\leq \eta_1 \}  \\
l'(\Phi\eta) &=& l'(\eta) + n-1- \#\{k>1| \eta_k\leq \eta_1 \} \\
a(s_i\eta) &=& a(\eta) \qquad l'(s_i\eta) = l'(\eta) \qquad
l(s_i\eta) = l(\eta) + 1 \quad\mbox{{\rm for} $\eta_i > \eta_{i+1}$} 
\end{eqnarray*}
\end{lemma}
A consequence of the first two equations in the final line is that
\begin{equation}\label{tu.10}
l'(\eta) = l'(\eta^+) = b(\eta^+), \qquad
a(\eta) = a(\eta^+) = b((\eta^+)')
\end{equation}
where $(\eta^+)'$ denotes the partition conjugate to $\eta^+$.
\subsection{Raising Operators and Lowering Operators}
There are two distinct raising operators which have a very simple action 
on non-symmetric Macdonald polynomials. Define \cite{knop96c,forr97b}
\begin{eqnarray}
\Phi_1 &:=& x_n\,\omega ,  \nonumber\\ \label{raise.1}
\Phi_2 &:=& x_n\,T_{n-1}^{-1}\cdots T_2^{-1} T_1^{-1}
\end{eqnarray}
A direct calculation reveals that for $i=1,2$
\begin{eqnarray*}
Y_n\,\Phi_i &=& q\Phi_i\,Y_1  \\
Y_j\,\Phi_i &=& \Phi_i\, Y_{j+1}\hspace{3cm} 1\leq j\leq n-1
\end{eqnarray*}
whence $\Phi_i E_{\eta}$ is a constant multiple of 
$E_{\Phi\eta}$, where $\Phi\eta:=(\eta_2,\ldots,\eta_n,\eta_1+1)$. This
constant is determined by looking at the coefficient of $x^{\Phi\eta}$
with the result that 
\begin{eqnarray*}
\Phi_1\,E_{\eta} &=& q^{\eta_1}\, E_{\Phi\eta}, \\
\Phi_2\,E_{\eta} &=& t^{-\#\{i|\eta_i\leq \eta_1\}}\, E_{\Phi\eta}
\end{eqnarray*}
{\it Remark.}\quad
These operators are simply related via $\Phi_1 = t^{n-1}\Phi_2\,Y_1$.
Of course any function of the operators $Y_i$ multiplied by $\Phi_1$
will be a raising operator for the non-symmetric Macdonald polynomials 
but these two are in some sense the simplest.

In a similar manner, one can use the $q$-Dunkl operators (\ref{dunk.old}) 
to construct lowering operators as follows,
\begin{eqnarray}
\Psi_1 &:=& \omega^{-1}\,D_n , \nonumber\\ 
\Psi_2 &:=& T_1 T_2\cdots T_{n-1}\,D_n \label{lower.1}
\end{eqnarray}
$\Psi_2$ was introduced previously in \cite{forr97b}.  
These operators intertwine with the Cherednik operators as 
\begin{eqnarray*}
Y_1\,\Psi_i &=& q^{-1}\Psi_i\,Y_n \\ 
Y_j\,\Psi_i &=& \Psi_i\,Y_{j-1}\hspace{3cm} 2\leq j\leq n
\end{eqnarray*}
and it is seen that
\begin{eqnarray*}
\Psi_1\,E_{\eta} &=& q^{-\eta_n+1}(1-t^{n-1+\bar{\eta}_n})\, 
E_{\Psi\eta}, \\
\Psi_2\,E_{\eta} &=& t^{\#\{i|\eta_i < \eta_n\}}
(1-t^{n-1+\bar{\eta}_n})
E_{\Psi\eta}
\end{eqnarray*}
where $\Psi\eta:=(\eta_n-1,\eta_1,\ldots,\eta_{n-1})$.

\subsection{Kernel}
Let\,\, $\widetilde{}$\,\, denote the involution on the ring of polynomials
with coefficients in $\cc(q,t)$, which acts on the the coefficients by
sending $q\mapsto q^{-1}$, $t\mapsto t^{-1}$, and extend it to act on
operators in the obvious way. Define the kernel 
\begin{equation}\label{calk}
{\cal K}_A(x;y;q,t) = \sum_{\eta} q^{a(\eta)} t^{(n-1)|\eta| - l'(\eta)}
\frac{d_{\eta}}{d'_{\eta}e_{\eta}} E_{\eta}(x)
\widetilde{E}_{\eta}(y)
\end{equation}
It follows from (\ref{qu}) that
this kernel is related to the previously introduced kernel 
\cite{forr97b}
\begin{equation}\label{kaa}
K_A(x;y;q,t) = \sum_{\eta} \frac{d_{\eta}}{d'_{\eta}e_{\eta}} 
E_{\eta}(x) \widetilde{E}_{\eta}(y)
\end{equation}
((\ref{kaa}) was denoted by ${\cal K}_A$ in \cite{forr97b}, but for the present
purpose it is desirable to use this notation for (\ref{calk})) 
by means of 
\begin{equation} \label{triste}
\widetilde{\cal K}_A(x;y;q,t) = K_A(-qy;x;q,t) 
\end{equation}

The kernel ${\cal K}_A(x;y;q,t)$ satisfies the following properties
\begin{thm} \label{kern.prop}
\begin{eqnarray*}
(a) &&(T_i^{\pm 1})^{(x)}\,{\cal K}_A(x;y;q,t) = \left(\widetilde{T_i^{\mp 1}}
\right)^{(y)}\,{\cal K}_A(x;y;q,t) \\
(b) &&\Psi_1^{(x)}\,{\cal K}_A(x;y;q,t) = \widetilde{\Phi_2}^{(y)}
\,{\cal K}_A(x;y;q,t) \\
(c) && \nd_i^{(x)}\,{\cal K}_A(x;y;q,t) = y_i \,{\cal K}_A(x;y;q,t)
\end{eqnarray*}
\end{thm}

\noindent
{\it Proof.} \quad
The proof of this result follows the same line of reason as in
\cite[Thm 5.2]{forr97b}, using the facts that
\begin{eqnarray}\label{i.1a}
x_i &=& t^{-n+i}\,\widetilde{T_i^{-1}}\cdots \widetilde{T_{n-1}^{-1}}
\widetilde{\Phi_2} \widetilde{T_1}\cdots \widetilde{T_{i-1}}  \\
\nd_i &=& t^{-n+i} T_{i-1}^{-1}\cdots T_1^{-1}\,\Psi_1\, 
T_{n-1}\cdots T_i \label{i.2a}
\end{eqnarray}
\hfill $\Box$

We recall from \cite{forr97b} that the analogue of property (c) for
the kernel
$K_A(x;y;q,t)$ is
\begin{equation}\label{c.old}
D_i^{(x)} K_A(x;y;q,t) = y_i  K_A(x;y;q,t).
\end{equation}
A feature of both property (c) and (\ref{c.old}) is that the $q$-Dunkl
operator
$\nd_i$ (resp.
$D_i$) act on the left set of variables {\it only}. However, by applying
the operation $\widetilde{\,\,}$ and using
(\ref{triste}), we can form similar identities where they act on 
the right set of variables, namely

\begin{cor}
\begin{eqnarray}
(\widetilde{\nd_i})^{(x)}\; K_A(z;x;q,t) &=&
-q^{-1}z_i\; K_A(z;x;q,t) \label{i.1}\\
(\widetilde{D_i})^{(y)} \;{\cal K}_A(x;y;q,t) &=&
-qx_i \;{\cal K}_A(x;y;q,t) \label{i.2}
\end{eqnarray}
\end{cor}
 
\subsection{First isomorphism}
Returning to the proof of Thm.~\ref{nata}, we claim that it follows
from the subsequent

\begin{prop} \label{pferd}
The map $\phi$ defined by
\begin{equation} \label{iso.1}
\phi(\widetilde{Y_i^{-1}}) = Y_i, \quad \phi(x_i)=e_i, \quad
\phi(\widetilde{T_i^{\pm 1}}) = T_i^{\mp 1} .
\end{equation}
is an  algebra isomorphism.
\end{prop}
{\it Proof.}\quad
As the relations between $\widetilde{Y_i^{-1}}$ and $x_j$ are somewhat
unsightly, it is more convenient to consider the operator
$\widetilde{\omega^{-1}}$ (from which the $\widetilde{Y_i^{-1}}$ can be
constructed) which possess simpler relations with the $x_j$
(c.f. (\ref{id.2})). Indeed by defining 
\begin{equation} \label{iso.2}
\phi(\widetilde{\omega^{-1}}) = T_1\cdots T_{n-1}\,\omega\,T_1^{-1}
\cdots T_{n-1}^{-1}
\end{equation}
we obtain the relation $\phi(\widetilde{Y_i^{-1}}) = Y_i$ as a consequence.

The explicit relations satisfied by the algebra $\{\widetilde{\omega^{-1}},
\widetilde{T_i^{\pm 1}}, x_i\}$ can be obtained by applying 
\,\, $\widetilde{}$\,\, to the relations (\ref{tdefs.1})--(\ref{tdefs.4}),
(\ref{id.2}). The fact that
the map $\phi$ as defined by (\ref{iso.1}), (\ref{iso.2}) is an isomorphism
then follows by a standard calculation. \hfill $\Box$
\vspace{.4cm}

\noindent
{\bf Proof of Thm.~\ref{nata}.}\quad We know that 
$E_{\eta}(x;q^{-1},t^{-1})$ is an
eigenfunction of $\widetilde{Y_i^{-1}}$. In principle this can be 
shown by utilizing the commutation relations amongst the operators
$\{\tilde{Y}_i^{-1}, \tilde{T}_i^{\pm 1}, x_i\}$ to move
the operators $\widetilde{Y_i^{-1}}$ through the terms in 
$E_{\eta}(x;q^{-1},t^{-1})$, until obtaining
$\widetilde{Y_i^{-1}} \cdot 1 = t^{-n+1} \cdot 1$. By adopting this viewpoint 
in  the eigenvalue equation (considered as an
operator equation)
$$
\widetilde{Y_i^{-1}}\,E_{\eta}(x;q^{-1},t^{-1})\,.\,1 = t^{\bar{\eta}_i}\,
E_{\eta}(x;q^{-1},t^{-1})\,.\,1
$$
and applying the map $\phi$ to both sides
it then follows from Lemma \ref{pferd} that 
$E_{\eta}(e;q^{-1},t^{-1}).1$ is an eigenfunction
of $\phi(\widetilde{Y_i^{-1}}) = Y_i$, with leading order term $x^{\eta}$
and hence must be proportional to $E_{\eta}(x;q,t)$. 

To determine the proportionality constant $\alpha_{\eta}(q,t)$ say,
it follows from
the action of $T_i$ given by (\ref{t.action}) that 
$$
e_1^{\eta_1} e_2^{\eta_2} \cdots e_n^{\eta_n}\,.\,1 = 
q^{f(\eta)}\,t^{g(\eta)}\, x^{\eta} + \sum_{\nu<\eta} b_{\eta\nu} x^{\nu}
$$
where $f(\eta) = \sum_i \bin{\eta_i}{2}$ and
\begin{eqnarray}
g(q) &=& \sum_{i=1}^n(\eta_i-1)(i-1) + \sum_{i=0}^{\eta_{n-1}} \chi(\eta_n\leq
i)  + \sum_{i=0}^{\eta_{n-2}} \chi(\eta_n\leq i) + \chi(\eta_{n-1}\leq i)
\nonumber \\
&&+ \cdots + \sum_{i=0}^{\eta_1} \chi(\eta_n\leq i) + \cdots
+ \chi(\eta_2 \leq i)
\end{eqnarray}
where $\chi(P)=1$ if $P$ is true, and zero otherwise. 
The simplification $g(q) = \sum_i (n-i)\eta_i^+ - \ell(w_{\eta})$ then 
follows from the above expression by induction on $\ell(w_{\eta})$. 
\hfill $\Box$

\setcounter{equation}{0}
\section{Al-Salam\&Carlitz polynomials}

The isomorphism $\phi$ introduced in the previous section can be
generalized to another isomorphism $\psi_a$ such that $\psi_a(x_i)$
includes not just degree-raising parts, but degree-preserving and lowering
parts as well. It will turn out that this isomorphism is precisely
what is needed to obtain non-symmetric analogues of the 
Al-Salam\&Carlitz polynomials in the same way as was done for the
Hermite case. 

As previously mentioned, the symmetric ASC polynomials $V_\lambda^{(a)}$
can be defined
via their orthogonality with respect to the inner product (\ref{innerv}).
We remark that under this inner product we have the important 
result that the adjoint	operators of $T_i^{\pm 1}$, $\omega$ are 
given by
\begin{equation} \label{adjoint}
(T_i^{\pm 1})^{\ast} = T_i^{\pm 1}, \qquad
(\omega^{-1})^{\ast} = \frac{t^{n-1}}{aq}\:\omega\,(x_1-q)(x_1-aq)
\end{equation}
The ASC polynomials  $V_\lambda^{(a)}$
can equivalently be defined by means of the 
generating function  \cite{forr97c}
$$
\prod_{i=1}^n {1 \over \rho_a(t^{-(n-1)}x_i;q)} {}_0{\psi}_0(x;y;q,t)
= \sum_{\lambda} {(-1)^{|\lambda|} q^{b(\lambda')}V^{(a)}_\lambda(y;q,t)
P_\lambda(x;q,t) \over d_\lambda'(q,t) P_\lambda(1,t,\dots,t^{n-1};q,t)}.
$$
Here, $\rho_a(x) := (x;q)_{\infty} (ax;q)_{\infty}$, $b(\lambda)$ is
defined by (\ref{def.n}) and 
\begin{eqnarray} 
P_\lambda(1,t,\dots,t^{n-1};q,t) &=& t^{l(\lambda)}\; \prod_{s\in\lambda}
\frac{(1-q^{a'(s)}t^{n-l'(s)})}{(1-q^{a(s)}t^{l(s)+1})} \nonumber\\
\label{kan.ker2}
{}_0{\psi}_0(x;y;q,t) &:=& \sum_{\lambda} {(-1)^{|\lambda|}
q^{b(\lambda)} \over
d_{\lambda}'(q,t)
P_{\lambda}(1,t,\dots,t^{n-1};q,t)} P_\lambda(x;q,t) P_\lambda(y;q,t).
\end{eqnarray}
This latter kernel was previously introduced by Kaneko \cite{kaneko93a}
in connection with hypergeometric solutions of systems of $q$-difference 
equations.

Similarly the  ASC polynomials  $U_\lambda^{(a)}$ can be defined by
the generating function \cite{forr97c}
\begin{equation}\label{gfmu}
\rho_a(x_1;q) \cdots \rho_a(x_n;q) \,
{}_0{\cal F}_0(x;y;q,t)
= \sum_{\kappa} {t^{b(\kappa)}U_\kappa^{(a)}(y;q,t)
P_\kappa(x;q,t) \over h_\kappa'(q,t) P_\kappa(1,t,\dots,t^{n-1};q,t)}
\end{equation}
where the hypergeometric function ${}_0{\cal F}_0$ is defined by
\begin{equation} \label{cano}
{}_0{\cal F}_0(x;y;q,t) := \sum_{\kappa} {t^{b(\kappa)} \over
d_\kappa'(q,t) P(1,t,\dots,t^{n-1};q,t)} P_\kappa(x;q,t) P_\kappa(y;q,t).
\end{equation}

\subsection{Second isomorphism}

Consider the involution \,\,$\hat{}$\,\, on polynomials and operators
defined in the statement of Corollary \ref{cor1.3}. 
The operator $E_i$ introduced in (\ref{e.def}) has its origins in this
involution, namely, 
\begin{equation}\label{def.E}
E_i:=(\hat{D}_{n+1-i})^{\ast}:= (-{1 \over q}{\nd}_{i})^{\ast}.
\end{equation} The form (\ref{e.def})
follows from (\ref{def.E}) by making use of the adjoint formulae 
(\ref{adjoint}). Regarding the algebra satisfied by the $E_i$ and
operators such as the $T_i$, note that
 application of 
the adjoint operation \,\,${}^{\ast}$\,\, to the 
relations involving ${\cal D}_i$, $T_i$ (\ref{nd.id})
gives, in place of the first relation in
(\ref{nd.id}) for example, 
\begin{equation}\label{tet}
T_i^{-1}\,E_i\,T_i^{-1} = t^{-1} E_{i+1}.
\end{equation}

Now consider the  mapping 
$
\psi_a : \{\widetilde{\omega}^{-1}, \widetilde{T_i}, x_i,
\widetilde{\nd}_i \} \longrightarrow \{ \omega, T_i, x_i, \nd_i \}
$
defined by
\begin{eqnarray}
\psi_a(x_i) &=& E_i, \nonumber\\
\psi_a(\widetilde{\omega}^{-1}) &=&
T_1\cdots T_{n-1} \left( Y_n + (1+a)t^{1-n} D_n + at^{2-2n}D_n Y_n D_n
\right) ,\nonumber \\
\psi_a(\widetilde{T}_i^{-1}) &=& T_i , \nonumber\\
\psi_a(\widetilde{\nd}_i) &=& -at^{n+1-2i} T_{i-1}\cdots T_1 T_1\cdots
T_{i-1}\,E_i^{\ast}\,T_i^{-1}\cdots T_{n-1}^{-1} T_{n-1}^{-1} \cdots
T_i^{-1} \label{iso.3}
\end{eqnarray}
 Then Theorem \ref{yanan} will follow from 

\begin{prop} \label{natillas}
The map $\psi_a$ is an algebra isomorphism.
\end{prop}

\noindent
{\it Proof.} \quad
The proof of this result consists of checking that the operators $\psi_a(u)$
given in (\ref{iso.3}) satisfy the same relations as the original operators
$u$, given by (\ref{tdefs.1})--(\ref{tdefs.4}), (\ref{id.2}) and
(\ref{nd.id}), (after application of the involution \,\,
$\widetilde{}$\,\, ). 
For example, the first formula in (\ref{2.1}), after application of the
involution \,\, $\widetilde{}$\,\, , reads
$$
\tilde{T}_i^{-1} x_{i+1} = t x_i \tilde{T}_i.
$$
Now applying the mapping  $\psi_a$ gives
$$
T_i E_{i+1} = t E_i T_i^{-1}.
$$
But this is equivalent to (\ref{tet}) so the algebra is indeed preserved.  
The calculations involved in checking the other relations are 
typically more tedious; however they are similar to those
undertaken in \cite{forr97b}, and so for brevity will be omitted.
\hfill $\Box$

\vspace{.2cm}
As with the relationship between Prop.~\ref{pferd}
and the proof of Thm.~\ref{nata} we are in a position to complete
the \\[2mm]
{\bf Proof of Thm.~\ref{yanan}}\quad From Thm~\ref{nata}, and the 
definition
(\ref{e.def}) of the operators $E_i$ it follows that  $E_{\eta}^{(V)}$ 
has leading term $E_{\eta}(x;q,t)$. In addition, it follows from
(\ref{iso.3}) that
\begin{equation}\label{3.4a}
\psi_a(\widetilde{{Y_i}^{-1}}) = Y_i + (1+a)t^{1-n}D_i + a t^{2-2n}
D_i\,Y_i^{-1}\,D_i
\end{equation}
and from Prop.~\ref{natillas}, that these are eigenoperators for the
non-symmetric ASC polynomials defined by (\ref{non.v}).
The corresponding eigenvalue is simply $t^{\bar{\eta}_i}$. 
By writing these operators out explicitly, it is seen that they are
self-adjoint w.r.t.~the inner product (\ref{innerv}). Hence by standard 
arguments, the polynomials (\ref{non.v}) are orthogonal w.r.t.~(\ref{innerv}).
\hfill $\Box$

\subsection{Normalization}
The images of the raising and lowering operators (\ref{raise.1}), 
(\ref{lower.1})
(after application of \,\,$\widetilde{}$\, ) under the map $\psi_a$ are 
guaranteed, by virtue of Prop.~\ref{natillas}, to be raising and lowering 
operators for the polynomials $E_{\eta}^{(V)}(x)$. 

In particular, using (\ref{i.2a}) and (\ref{iso.3}) we see that 
$$
\psi_a(\widetilde{\Psi_1}) = aq^{-1}t^{1-n}\,\Psi_1
$$
so that $\Psi_1$ remains a raising operator for the polynomials
$E_{\eta}^{(V)}$. By examination of the leading terms, we must have
\begin{equation} \label{sacerdote.1}
\Psi_1\, E_{\eta}^{(V)} = q^{\eta_n+1} \frac{d'_{\eta}}{d'_{\Psi\eta}}\;
E_{\Psi\eta}^{(V)}
\end{equation}
Also, use of (\ref{i.1a}) and (\ref{iso.3}) gives 
$$
\psi_a(\widetilde{\Phi_2}) = -q^{-1}\,\Psi_1^{\ast}
$$
so that $\Psi_1^{\ast}$ is a raising operator for $E_{\eta}^{(V)}$.
Indeed,
\begin{equation} \label{sacerdote.2}
\Psi_1^{\ast}\,E_{\eta}^{(V)} = a^{-1} t^{n-1} q^{\eta_1+1} \, 
E^{(V)}_{\Phi\eta}
\end{equation}
By an argument similar to that used in \cite[Prop. 3.6]{forr96d}
it follows from (\ref{sacerdote.1}) and (\ref{sacerdote.2}) that
\begin{equation} \label{recurrence.1}
\inn{E_{\Phi\eta}^{(V)}}{E_{\Phi\eta}^{(V)}}^{(V)} = a t^{1-n} q^{-2\eta_1-1}\,
\frac{d'_{\Phi\eta}}{d'_{\eta}}\; \inn{E_{\eta}^{(V)}}{E_{\eta}^{(V)}}^{(V)}
\end{equation}
Also, we have
\begin{equation} \label{recurrence.2}
\inn{E_{s_i\eta}^{(V)}}{E_{s_i\eta}^{(V)}}^{(V)}
= \frac{(1-t^{\delta_{i\eta}-1})
(1-t^{\delta_{i\eta}+1})}{t(1-t^{\delta_{i\eta}})^2} \;
\inn{E_{\eta}^{(V)}}{E_{\eta}^{(V)}}^{(V)}
\end{equation}
The solution of the recurrence relations (\ref{recurrence.1}), 
(\ref{recurrence.2}) gives
\begin{prop}
\begin{equation}\label{norm.v}
{\cal N}^{(V)}_{\eta} := \inn{E_{\eta}^{(V)}}{E_{\eta}^{(V)}}^{(V)}
= \left(aq^{-1}t^{2-2n}\right)^{|\eta|}\,q^{-2a(\eta)} t^{l(\eta)+l'(\eta)}
\frac{d'_{\eta}\,e_{\eta}}{d_{\eta}}\;{\cal N}^{(V)}_0
\end{equation}
where for $t=q^k$, \cite{forr97c}
$$
{\cal N}^{(V)}_0 = (1-q)^n a^{kn(n-1)/2}
t^{-2 k \left ({n \atop 3} \right ) -k
\left ({n \atop 2} \right )}
\prod_{l=1}^n {(q;q)_{kl} \over (q;q)_k}.
$$
\end{prop}

By using the formulas (\ref{uv}), (\ref{qiuv}) and
(\ref{deltain}) we see that the norm ${\cal
N}^{(U)}_{\eta}$ of the non-symmetric ASC polynomials $E_\eta^{(U)}$ with
respect to the inner product (\ref{inneru}) is given by simply replacing
$q,t$ by $q^{-1}, t^{-1}$ in ({\ref{norm.v}). Use of (\ref{qu}) then gives
\begin{cor}
\begin{equation}\label{norm.u}
{\cal N}^{(U)}_{\eta} := \inn{E_{\eta}^{(U)}}{E_{\eta}^{(U)}}^{(U)}
= \left(at^{n-1}\right)^{|\eta|}\,q^{a(\eta)} t^{-l(\eta)}
\frac{d'_{\eta}\,e_{\eta}}{d_{\eta}}\;{\cal N}^{(U)}_0
\end{equation}
where for $t=q^k$, \cite{forr97c}
$$
{\cal N}^{(U)}_0 = (1-q)^n (-a)^{kn(n-1)/2}
t^{ k \left ({n \atop 3} \right ) -{k-1 \over 2}
\left ({n \atop 2} \right )}
\prod_{l=1}^n {(q;q)_{kl} \over (q;q)_k}.
$$
\end{cor}

\subsection{Generating function}

The raising operator expression (\ref{non.v}) facilitates the 
derivation of the generating function for the non-symmetric
ASC polynomials. 
Also required will be the $q$-symmetrization of (\ref{calk}).
\begin{prop}
Let \cite{mac95} $U^+ = \sum_{\sigma} T_{\sigma}$ where
$T_{\sigma}:=T_{i_1}\cdots T_{i_p}$ for a reduced word decomposition
$\sigma=s_{i_1}\cdots s_{i_p}$. We have
\begin{equation} \label{kaneko.kernel}
(U^+)^{(x)}\; {\cal K}_A(x;y;q,t) =  [n]_t! {}_0\psi_0(x;-t^{n-1}y;q,t)
\end{equation}
where ${}_0\psi_0$ is defined by (\ref{kan.ker2}).
\end{prop}

\noindent
{\it Proof.} \quad We remark that this is the analogue of the result
\cite[Prop.~5.4]{forr97b}
\begin{equation}\label{ukf}
(U^+)^{(x)}\;  K_A(x;y;q,t) =  [n]_t! {}_0F_0(x;y;q,t)
\end{equation}
In fact in our proof of (\ref{kaneko.kernel}) we will use the formula
\begin{equation}
U^+\,E_{\eta}(x) = [n]_t! t^{l(\eta)}\,\frac{e_{\eta}}{P_{\lambda}
(t^{\delta})  d_{\eta}} \, P_{\lambda}(x), \qquad \lambda=\eta^+
\label{sym.2}
\end{equation}
which was deduced \cite[eqs.~(5.8)\&(5.18)]{forr97b} 
as a corollary of (\ref{ukf}). 
Thus we apply $U^+$ to (\ref{calk}) and use (\ref{sym.2}) to compute its
action. Simplifying the result using the first equation in (\ref{tu.10})
and the formula \cite{mac95}
\begin{equation}
P_{\lambda}(y)  =  \sum_{\eta:\eta^+=\lambda} \frac{d'_{\lambda}}
{d'_{\eta}}\: E_{\eta}(y), \label{sym.1}
\end{equation}
the result then follows.
\hfill$\Box$

\vspace{.2cm}
Consider now the generating function
$$
F_1(y ; z) = \sum_{\nu} A_{\nu} E^{(V)}_{\nu}(y)\, \widetilde{E}_{\nu}(z)
$$
where 
\begin{equation}\label{A}
A_{\nu} = (a/q)^{|\nu|}\,\frac{{\cal N}_0^{(V)}}{\alpha_{\nu}(q,t)
{\cal N}_{\nu}^{(V)}}
= q^{a(\nu)} t^{(n-1)|\nu| - l'(\nu)}\:
\frac{d_{\nu}}{d'_{\nu}e_{\nu}}
\end{equation}
Here we have used the fact that $l(\eta) = l(\eta^+) + \ell(w_{\eta})$
to rewrite $\alpha_{\eta}(q,t)$ as defined by (\ref{cuerda}) as
$$
\alpha_{\eta}(q,t) = q^{a(\eta)}\,t^{(n-1)|\eta| - l(\eta)} .
$$
Clearly
$$
\inn{F_1(y;z)}{E^{(V)}_{\eta}(y)}_y^{(V)} = (a/q)^{|\eta|}\,
\frac{{\cal N}_0^{(V)}}
{\alpha_{\eta}(q,t)}\; \widetilde{E}_{\eta}(z)
$$
Next note the integration formula 
\begin{eqnarray}
\inn{{\cal K}_A(y;z)}{1}_y^{(V)}
&=& \frac{1}{[n]_t!}\,\inn{U_y^+\,{\cal K}_A(y;z)}
{1}_y^{(V)}  \nonumber\\
&=& \inn{{}_0\psi_0(y;-t^{n-1}z)}{1}_y^{(V)} = {\cal N}_0^{(V)} \prod_{i=1}^n
\rho_a(-z_i) \label{seso}
\end{eqnarray} 
which follows from the symmetrization formula (\ref{kaneko.kernel}),
the fact that $U^+_y$ is self adjoint w.r.t.~$\langle \, , \,
\rangle_y^{(V)}$
and an integral formula for the kernel ${}_0\psi_0(y;z)$ given in
\cite[Prop 4.8]{forr97c}, and
consider the generating function
$$
F_2(y ; z) = \prod_{i=1}^n \frac{1}{\rho_a(-z_i)}\; {\cal K}(y;z)
$$
We have 
\begin{eqnarray*}
\inn{F_2(y;z)}{E_{\eta}^{(V)}(y)}_y^{(V)} 
&=& \frac{(-a)^{|\eta|}}{\alpha_{\eta}(q,t)}\prod_i\frac{1}{\rho_a(-z_i)}
\inn{{\cal K}(y;z)}{\widetilde{E}_{\eta}(E^{(y)})}_y^{(V)} \\
&=& \frac{(a/q)^{|\eta|}}{\alpha_{\eta}(q,t)}\prod_i\frac{1}{\rho_a(-z_i)}
\inn{\widetilde{E}_{\eta}(\nd^{(y)})\,{\cal K}(y;z)}{1}_y^{(V)}  \\
&=& \frac{(a/q)^{|\eta|}}{\alpha_{\eta}(q,t)} \prod_i\frac{1}{\rho_a(-z_i)}
\widetilde{E}_{\eta}(z)\, \inn{{\cal K}(y;z)}{1}_y^{(V)} \\
&=& (a/q)^{|\eta|} \frac{{\cal N}_0^{(V)}}{\alpha_{\eta}(q,t)}
\widetilde{E}_{\eta}(z)
\end{eqnarray*}
In the above chain of equalities, we have used (\ref{non.v}), 
(\ref{def.E}), the kernel property Thm.~\ref{kern.prop} (c)
and (\ref{seso}) respectively. The non-symmetric ASC polynomials 
$E_{\eta}^{(V)}(y)$ are a complete basis for polynomials in $y$ and
hence from above we have $F_1=F_2$. That is, we have the generating
function for non-symmetric ASC polynomials $ E^{(V)}_{\nu}$.
\begin{prop} With $A_\nu$ given by (\ref{A})
\begin{equation} \label{gen.fn}
\prod_{i=1}^n \frac{1}{\rho_a(-z_i)}\; {\cal K}_A(y;z) =
\sum_{\nu} A_{\nu} E^{(V)}_{\nu}(y)\, \widetilde{E}_{\nu}(z) .
\end{equation}
\end{prop}

We remark that this generating function could also be derived 
in a manner similar to that done in the symmetric case \cite{forr97c},
namely by applying the operator $(\widetilde{Y_i^{-1}})^{(z)}$ to
both sides of (\ref{gen.fn}) and deducing that $E^{(V)}_{\eta}(y)$
is an eigenfunction of
\begin{equation}
h_i = \psi_a(\widetilde{Y_i^{-1}}) = Y_i\,T_{i-1}\cdots T_1 \left(1+\nd_1 
\right) \left(1+a\nd_1 \right) T_1^{-1}\cdots T_{i-1}^{-1}
\end{equation}
with leading term $E_{\eta}(y)$ (some manipulation using (\ref{nd.id})
and (\ref{dunk.new}) cast this into the form given in (\ref{def.hop})). 
Note also that by applying the operation $\hat{\,\,}$ with the
respect to the $y$-variables in (\ref{gen.fn}) and using the formula
(\ref{triste}) as well as
$$
{1 \over \rho_a(x;q)} \bigg |_{q \mapsto q^{-1}} = \rho_a(qx;q),
$$
(see e.g.~\cite{forr97c})
we deduce the generating function formula for the polynomials
$ E^{(U)}_{\nu}$
\begin{cor}
\begin{equation}\label{genU}
\prod_{i=1}^n \rho_a(z_i) \, K_A(z;y^R;q,t) =
\sum_\nu  {d_\nu \over d_\nu' e_\nu}
E_\nu^{(U)}(y) E_\nu(z)
\end{equation}
\end{cor}

The generating function formulas in turn imply 
a further class of operator formulas
relating the ASC polynomials and the non-symmetric Jack polynomials
(c.f.~\cite[eqs.~(3.9)\&(3.10)]{forr97c}).
\begin{cor}
We have
\begin{eqnarray}
E_{\eta}^{(V)}(y)&  = & \prod_{i=1}^n \frac{1}{\rho_a(-\nd^{(y)}_i)}\;
E_{\eta}(y) \label{b1}\\
E_{\eta}^{(U)}(y)&  = & \prod_{i=1}^n \rho_a\Big (- q 
\widetilde{{\nd}^{(y)}}_i\Big 
)\;
\widetilde{E}_{\eta}(y^R) \label{b2}.
\end{eqnarray}
\end{cor}

\noindent
{\it Proof.} \quad 
The first identity follows
from (\ref{gen.fn}) by using Thm.~\ref{kern.prop} (c)
and comparing coefficients of $\widetilde{E}_{\eta}(z)$, while the
second identity follows similarly from (\ref{genU}) and (\ref{i.1}).
\hfill$\Box$

\vspace{.2cm}
As further applications of the generating functions we will present
some evaluation formulas for $E_\eta^{(V)}$ at the special points
$t^{\bar{\delta}-n + 1}$ and $at^{\bar{\delta}-n + 1}$, where
$t^{\bar{\delta}} := (1,t,t^2,\dots,t^{n-1})$.

\begin{prop}
We have
\begin{eqnarray}
E_\eta^{(V)}(t^{\bar{\delta}-n + 1})&  = & (-a)^{|\eta|} q^{-a(\eta)}
t^{l'(\eta) - (n-1) |\eta|} E_\eta(t^{\bar{\delta}}) \label{n.1}\\
E_\eta^{(V)}(at^{\bar{\delta}-n + 1})&  = & (-1)^{|\eta|} q^{-a(\eta)}
t^{l'(\eta) - (n-1) |\eta|} E_\eta(t^{\bar{\delta}}) \label{n.2}
\end{eqnarray}
where
\begin{equation}\label{n.3}
E_\eta(t^{\bar{\delta}}) = t^{l(\eta)} {e_\eta \over d_\eta}.
\end{equation}
\end{prop}

\noindent
{\it Proof.} \quad The formula (\ref{n.3}) is a special case of a result
of Cherednik \cite{cher95b} (see also \cite{noumi96c}). For the derivation 
of (\ref{n.1}) and (\ref{n.2}) we follow the strategy of the proof of the
analogous result in the symmetric case \cite[Prop.~4.3]{forr97c}. First,
note from the definition (\ref{maldonado}) that in general
$$
T_i f(t^{\bar{\delta}}) = t  f(t^{\bar{\delta}}),
$$
and so
$$
U^+  f(t^{\bar{\delta}}) = (U^+ \, 1)  f(t^{\bar{\delta}}) =
[n]_t!  f(t^{\bar{\delta}}).
$$
Use of this latter formula in (\ref{kaneko.kernel}) with $y=t^{\bar{\delta}}$
gives
\begin{equation}\label{n.4}
{\cal K}_A (t^{\bar{\delta}};z;q,t) = {}_0 \psi_0(t^{\bar{\delta}};
-t^{n-1}z;q,t) = \prod_{i=1}^n (-t^{n-1}z_i;q)_\infty,
\end{equation}
and similarly, from (\ref{ukf})
\begin{equation}\label{n.5}
K_A (t^{\bar{\delta}};z;q,t) = {}_0 F_0(t^{\bar{\delta}};z;q,t) =
 {1 \over \prod_{i=1}^n (z_i;q)_\infty},
 \end{equation}
where the final equalities in (\ref{n.4}) and (\ref{n.5}) are known
results \cite{macunp1, kaneko96a}. Now set $y=t^{\bar{\delta}-n+1}$ in the 
generating function  (\ref{kaneko.kernel}). Use of (\ref{n.4}) with
$z$ replaced by $t^{-n+1}z$, and then use of (\ref{n.5}) allows the
l.h.s.~of the resulting expression to be written
$$
{1 \over \prod_{i=1}^n  (-az_i;q)_\infty} = K_A(t^{\bar{\delta}};
-az;q,t) = \sum_{\eta} {(-a)^{|\eta|} d_\eta \over d_\eta' e_\eta}
E_\eta(t^{\bar{\delta}}) \widetilde{E}_\eta(z).
$$
Equating with $ \widetilde{E}_\eta(z)$ on the r.h.s.~of the resulting
expression gives (\ref{n.1}). The formula (\ref{n.2}) follows similarly,
by substituting $y = at^{\bar{\delta}-n+1}$ in  (\ref{kaneko.kernel}).
\hfill$\Box$

\subsection{Relationship to the symmetric ASC polynomials}
The non-symmetric ASC polynomials are related to the corresponding
symmetric ASC polynomials in an analogous way to the relationship
(\ref{sym.2}) between the non-symmetric and symmetric Macdonald 
polynomials.
\begin{prop}
Let 
$$
a_\eta(q,t) = [n]_t! t^{\ell(\eta)} {e_\eta \over
P_{\eta^+}(t^{\bar{\delta}}) d_\eta}.
$$
We have
\begin{eqnarray}
U^+ E_\eta^{(V)}(y)&  = & a_\eta(q,t) V_{\eta^+}^{(a)}(y;q,t) \label{s.1}\\
U^+ E_\eta^{(U)}(y)& = & a_\eta(q,t) U_{\eta^+}^{(a)}(y;q,t).
\end{eqnarray}
\end{prop}

\noindent
{\it Proof.} \quad Consider the action of the $U^+$ operator on
(\ref{b1}) and (\ref{b2}). From the first three equations of
(\ref{nd.id}) one can check that $T_i$ commutes with any symmetric
function of the ${\cal D}_i$. Thus the action of $U^+$ can be
commuted to act to the right of $\prod_i\rho_a(-{1 \over q} \tilde{\cal
D}_i)$ and $1/\prod_i\rho_a(-{\cal D}_i)$. Use of (\ref{sym.2})
then gives
\begin{eqnarray*}
U^+ E_\eta^{(V)}(y) &  = &  a_\eta(q,t) { 1 \over \prod_i\rho_a(-{\cal D}_i)}
P_{\eta^+}(y)  =   a_\eta(q,t) {  1 \over \prod_i\rho_a(q 
\widetilde{D}_i)}
P_{\eta^+}(y)
\\
U^+ E_\eta^{(U)}(y) &  = &  a_\eta(q,t) \prod_i\rho_a(- q
\widetilde{\cal D}_i) \, P_{\eta^+}(y)  =   a_\eta(q,t)
 \prod_i \rho_a(D_i) \, P_{\eta^+}(y),
\end{eqnarray*}
where in obtaining the first equality in the second formula we have
used the fact that $\widetilde{P}_\eta(y^R) = P_\eta(y)$, while the
second equalities in both formulas make use of (\ref{dunk.new})
and the fact that $P_{\eta^+}$ is symmetrical. But the resulting
operator formulas are precisely representations obtained in
\cite[eq.~(3.9)\&(3.10)]{forr97c} for the symmetric ASC polynomials.
\hfill$\Box$

\vspace{.2cm}
We can also relate the eigenoperators $h_i$ for the non-symmetric ASC
polynomials $E^{(V)}_\eta$ to the eigenoperator \cite[eq.~(3.28)]{forr97c}
\begin{equation} \label{form.2}
{\cal H} = t^{1-n}\sum_{i=1}^n Y_i^{-1} - (1+a) \sum_{i=1}^n t^{1-i}
D_i\,Y_i^{-1} + a \sum_{i=1}^n t^{1-i} D_i^2\,Y_i^{-1} +
a(1-t^{-1})\sum_{1\leq i<j\leq n} t^{1-i} D_j\,D_i\,Y_i^{-1}
\end{equation}
for the symmetric ASC polynomials $U_\lambda^{(a)}$.

\begin{prop}
Let $h_i$ be given by (\ref{def.hop}) and $\cal H$ by (\ref{form.2}).
When acting on symmetric functions
$$
\sum_{i=1}^n h_i = t^{1-n} \widetilde{\cal H}.
$$
\end{prop}

\noindent
{\it Proof.} \quad {}From Theorem \ref{yanan}, by summing over $i$ in
(\ref{def.hop}) we have
$$
\sum_{i=1}^n h_i \, E^{(V)}(x;q,t) = t^{1-n} e(\eta^+)  E^{(V)}(x;q,t),
$$
where $e(\eta^+) = \sum_{i=1}^n t^{\bar{\eta}_i} = 
\sum_{i=1}^n q^{\eta_i^+} t^{n-i}$. We would next like to apply the
operator $U^+$ to both sides of this eigenvalue equation. For this purpose
we require the fact that $T_i$ commutes with $\sum_{i=1}^n h_i$
(this follows from (\ref{ty}), and the fact that these same equations
apply with the $Y_i$ replaced by $D_i$). Thus, making use of
(\ref{s.1}), this operation gives
$$
\sum_{i=1}^n h_i \, V_{\eta^+}^{(a)}(x;q,t) = t^{1-n}  e(\eta^+)
 V_{\eta^+}^{(a)}(x;q,t).
$$
But from \cite{forr97c} we know that this same eigenvalue equation applies
with $\sum_{i=1}^n h_i$ replaced by $t^{1-n} \tilde{\cal H}$. The result
now follows from the fact that $\{V_{\eta^+}^{(a)} \}$ are a basis for
symmetric functions.

We remark that an alternative proof is to establish directly that when
acting on symmetric functions
\begin{eqnarray}
\sum_{i=1}^n \tilde{Y}_i^{-1} & = & \sum_{i=1}^n Y_i \\
 - \sum_{i=1}^n t^{-1+i} \tilde{D}_i \tilde{Y}_i^{-1}&  = & \sum_{i=1}^n D_i \\
\sum_{i=1}^n t^{-1+i} \tilde{D}_i^2 \tilde{Y}_i^{-1} + (1-t)
\sum_{1 \le i < j \le n} t^{-1+i}  \tilde{D}_j  \tilde{D}_i
\tilde{Y}_i^{-1} & = &
t^{1-n} \sum_{i=1}^n D_i Y_i^{-1} D_i.
\end{eqnarray}
\hfill$\Box$

\subsection{Non-symmetric shifted Macdonald polynomials}

In \cite{forr97c} it was observed that the symmetric ASC polynomials
$V^{(a)}_{\lambda}(x)$ coincide (up to a factor and change of
variables) with the shifted Macdonald polynomials when $a=0$. We show now 
that this behaviour carries over to the non-symmetric case.

Following Knop and Sahi \cite{knop96a,knop96d,sahi97a},
the non-symmetric shifted Macdonald polynomials $G_{\eta}(z)$ are defined
as the unique polynomial with expansion
$$
G_{\eta}(z;q,t) = \widetilde{E}_\eta(z) + \sum_{|\nu| < |\eta|} b_{\eta\nu}
\widetilde{E}_{\nu}(z)
$$
(in \cite{sahi97a} what we denote 
$ \widetilde{E}_\eta(z)$ is called the non-symmetric Jack polynomial), 
which vanishes at the points $z= t^{\bar{\xi}}$ for all compositions $\xi$
such that $|\xi| \leq |\eta|$. Here $t^{\bar{\xi}}$ is given by 
(\ref{e-val}). Equivalently \cite{knop96a,okoun96b} they can be
defined as eigenfunctions of the ``inhomogeneous'' Cherednik operators
$$
\Xi_i = \widetilde{Y}_i + \widetilde{D}_i  
$$
where the operators are defined with the variables $z_i$.
For such polynomials, Knop and Sahi defined a raising operator 
$\Phi_{KS}= (z_n - t^{1-n})\omega^{-1}$ with a simple action on 
$G_{\eta}(z;q,t)$. It is easily seen that
$$
\lim_{a\rightarrow 0} \frac{-a}{q} \Psi_1^{\ast} =  \widetilde{\Phi}_{KS}
\Bigg|_{z_i=t^{n-1}x_i}, \qquad
\lim_{a\rightarrow 0} h_i =  \widetilde{\Xi_i}
\Bigg|_{z_i=t^{n-1}x_i }
$$
which immediately implies the sought relationship between $G_\eta$ and
$E_\eta^{(V)}$.
\begin{prop}\label{EG0}
\begin{equation}\label{EG}
\left. E^{(V)}_{\eta}(x;q,t)\right|_{a=0} = t^{-(n-1)|\eta|}
G_{\eta}(t^{n-1}x;q^{-1},
t^{-1})
\end{equation}
or equivalently
\begin{equation}\label{EG1}
\left. E^{(U)}_{\eta}(x;q,t)\right|_{a=0} = t^{(n-1)|\eta|}
G_{\eta}(t^{1-n}x;q,
t)
\end{equation}

\end{prop}

One immediate application of (\ref{EG}) is the evaluation of
$G_\eta(0;q,t)$, which follows from (\ref{n.2}). This
is a special case of a result of Sahi \cite[Th.~1.1]{sahi97a}, in which
an evaluation formula is given for $G(\alpha t^{\bar{\delta}};q,t)$,
for a general scalar $\alpha$. In fact use of (\ref{EG}) also allows this
more general
evaluation formula to be deduced. 

\begin{prop}
With  $(\alpha)_\lambda^{(q,t)} := \prod_{s \in \lambda}
(t^{l'(s)} - q^{a'(s)} \alpha)$ we have
$$
G_\eta(t^{-\bar{\delta}}\alpha;q,t) = \alpha^{|\eta|}
(1/\alpha)_{\eta^+}^{(q,t)} t^{-(n-1)|\eta|} {e_\eta \over d_\eta}.
$$
\end{prop}

\noindent
{\it Proof.} \quad Choosing $a=0$ and $y=t^{n-1-\bar{\delta}}\alpha$ in 
(\ref{genU}), and using (\ref{n.5}) and (\ref{EG1}), we see that
$$
\sum_{\eta} \alpha^{-|\eta|} t^{(n-1)|\eta|} {d_\eta \over d_\eta' e_\eta}
G_\eta(t^{-\bar{\delta}}\alpha;q,t) E_\eta(z) =
\prod_{i=1}^n { (z_i/\alpha;q)_\infty \over (z_i;q)_\infty}.
$$
But we know that \cite{macunp1,kaneko96a}
$$
\prod_{i=1}^n { (z_i/\alpha;q)_\infty \over (z_i;q)_\infty} =
\sum_\lambda {(1/\alpha)_\lambda^{(q,t)} \over d_\kappa'} P_\kappa(z;q,t)
= \sum_\eta {(1/\alpha)_{\eta^+}^{(q,t)} \over d_\eta'} E_\eta(z).
$$
The result follows by equating coefficients of $E_\eta(z)$.
\hfill$\Box$

\subsection{$q$-binomial coefficients}
Sahi \cite{sahi97a} uses the polynomials $G_\eta$ to introduce
non-symmetric $q$-binomial coefficients
 $\left [ {\eta \atop \nu} \right ]_{q,t}$ according to
\begin{equation}
\left [ {\eta \atop \nu} \right ]_{q,t} := {G_\nu(t^{\bar{\eta}}) \over
G_\nu(t^{\bar{\nu}})}
\end{equation}
($\bar{\eta}_i$ is defined by (\ref{e-val})). Our generating function
characterization of the ASC polynomials, and thus by Proposition \ref{EG0}
of the polynomials $G_\eta$, makes it natural to extend Lassalle's
\cite{lass97a} definition of the symmetric $q$-binomial coefficients to the
non-symmetric case by defining the non-symmetric  $q$-binomial coefficients
$\left ( {\eta \atop \nu} \right )_{q,t}$ according to the generating
function formula
\begin{equation}\label{nsqb}
\widetilde{E}_\nu(x) \prod_{i=1}^n {1 \over (x_i;q)_\infty} =
\sum_\eta \left ( {\eta \atop \nu} \right )_{\! q,t}
t^{l(\eta) - l(\nu)} {d_\nu' \over d_\eta'} \widetilde{E}_\eta(x)
\end{equation}
We can then use the generating function (\ref{kaneko.kernel}) to
relate these binomial coefficients to the polynomials $G_\eta$.

\begin{prop}
With $\left ( {\eta \atop \nu} \right )_{q,t}$ defined by (\ref{nsqb}), we
have
\begin{equation}\label{tu.1}
{G_\eta(x) \over G_\nu(0)} = \sum_\nu \left ( {\eta \atop \nu} \right
)_{\! q^{-1},t^{-1}} {\widetilde{E}_\nu(x) \over G_\nu(0) }
\end{equation}
\end{prop}

\noindent
{\it Proof.} \quad Multiply both sides of (\ref{nsqb}) by 
$q^{a(\nu)} t^{(n-1)|\nu|-l'(\nu)}{d_\nu \over d_\nu' e_\nu} E_\nu(y)$
and sum over $\nu$, rewriting the l.h.s.~according to (\ref{kaneko.kernel}).
Now equate coefficients of  $\widetilde{E}_\nu(x)$ on both sides. The
result then follows upon using (\ref{n.3}) and  (\ref{EG}).
\hfill$\Box$

\vspace{.2cm}
Since (\ref{tu.1}) is a formula satisfied by the non-symmetric
$q$-binomial coefficients of Sahi \cite[Cor.~1.3]{sahi97a}, and this formula
suffices to implicitly define these coefficients, we have that
\begin{equation}
\left ( {\eta \atop \nu} \right )_{q,t} =
\left [ {\eta \atop \nu} \right ]_{q,t}.
\end{equation}

Finally, let us present some formulas relating the cefficients
$\left ( {\eta \atop \nu} \right )_{q,t}$ to their symmetric
counterparts $\left ( {\kappa \atop \mu} \right )_{q,t}$, which can be
characterized by either of the formulas \cite{lass97a,okoun97a}
\begin{eqnarray}\label{las}
P_{\mu}(x;q,t)\;\prod_{i=1}^n \frac{1}{(x_i;q)_{\infty}}&  = &
\sum_{\lambda}\left ({\lambda \atop \mu}\right)_{\!q,t}t^{b(\lambda)-b(\mu)}
\frac{d'_{\mu}}{d'_{\lambda}}\;P_{\lambda}(x;q,t). \\
\label{okeq}
{P_\lambda^*(y;q^{-1},t^{-1}) \over P_\lambda^*(0;q^{-1},t^{-1})}
& = &\sum_{\mu} \left ( {\lambda \atop \mu} \right )_{\!q,t}
 {P_\mu(yt^{\bar{\delta}};q,t) \over P_\lambda^*(0;q^{-1},t^{-1})}
\end{eqnarray}
Here $P_\lambda^*$ is the shifted Macdonald polynomial, which is related
to the symmetric ASC polynomial $V_\lambda^{(0)}$ by
\cite[Prop.~4.4]{forr97c}
\begin{equation}\label{psv}
P_\lambda^*(yt^{-\bar{\delta}+n-1};q^{-1},t^{-1}) =
 t^{(n-1)|\lambda|}  V_\lambda^{(0)}(y;q,t).
 \end{equation}
\begin{prop}
With $\eta^+ = \kappa$, $\nu^+ = \mu$,
\begin{eqnarray}\label{tu.5}
\sum_{\nu : \nu^+ = \mu} 
\left ( {\eta \atop \nu} \right )_{q,t} &  = & \left ( {\kappa \atop \mu} \right
)_{q,t} \\
{ d_\kappa' \over d_\mu'}
{P_{\kappa}(t^{\bar{\delta}}) \over P_{\mu}(t^{\bar{\delta}})}
{d_\nu' \over E_\nu(t^{\bar{\delta}})} \sum_{\eta:\eta^+ = \kappa}
\left ( {\eta \atop \nu} \right )_{q,t} { E_\eta(t^{\bar{\delta}})
\over d_\eta' } & = &  \left ( {\kappa \atop \mu} \right)_{q,t}
\label{tu.6}
\end{eqnarray}
\end{prop}

\noindent
{\it Proof.} \quad The proof follows the strategy given in \cite{forr97c}
for the proof of the corresponding results in the $q=t^\alpha, \, q \to 1$
limit (binomial coefficients associated with non-symmetric Jack
polynomials). For (\ref{tu.5}) we apply the $U^+$ operator to
(\ref{tu.1}), making use of (\ref{sym.2}) and (\ref{s.1}). Use of the
fact that
$$
{a_\nu \over E_\nu^{(V)}(0)} = {[n]_t! \over V_{\eta^+}^{(0)}(0;q,t)}
$$
and (\ref{psv}) then gives
$$
{ P_\lambda^*(xt^{-\bar{\delta}};q^{-1},t^{-1}) \over
 P_\lambda^*(0;q^{-1},t^{-1}) } =
 \sum_{\nu} \left ( {\eta \atop \nu} \right )_{q,t}
 {P_{\nu^+}(x;q,t) \over P_{\nu^+}^*(0;q^{-1},t^{-1})}.
$$
Comparison with (\ref{okeq}) implies (\ref{tu.5})
The identity (\ref{tu.6}) follows similarly, by applying $U^+$ to
(\ref{nsqb}) and comparing with (\ref{las}).
\hfill$\Box$

\vspace{5mm}\noindent
{\bf\Large Acknowledgements}\\[2mm]
The authors acknowledge the financial support of the Australian Research
Council.

\bibliographystyle{plain}

\begin{thebibliography}{10}

\bibitem{forr97c}
T.~H. Baker and P.~J. Forrester.
\newblock Multivariable {Al-Salam \& Carlitz} polynomials associated with type
  {$A$} $q$-{Dunkl} operators.
\newblock q-alg/9706006.

\bibitem{forr96d}
T.~H. Baker and P.~J. Forrester.
\newblock Non--symmetric {Jack} polynomials and integral kernels.
\newblock q-alg/9612003, to appear in {\it Duke J. Math.}

\bibitem{forr96c}
T.~H. Baker and P.~J. Forrester.
\newblock The {Calogero--Sutherland} model and polynomials with prescribed
  symmetry.
\newblock {\em Nucl. Phys.}, B 492:682--716, 1997.

\bibitem{forr97b}
T.~H. Baker and P.~J. Forrester.
\newblock A $q$-analogue of the type {$A$} {Dunkl} operator and integral
  kernel.
\newblock {\em Int. Math. Res. Not.}, 14:667--686, 1997.

\bibitem{cher91a}
I.~Cherednik.
\newblock A unification of the {Knizhnik--Zamolodchikov} and {Dunkl} operators
  via affine {Hecke} algebras.
\newblock {\em Inv. Math.}, 106:411--432, 1991.

\bibitem{cher94a}
I.~Cherednik.
\newblock Integration of quantum many--body problems by affine
  {Knizhnik--Zamolodchikov} equations.
\newblock {\em Commun. Math. Phys.}, 106:65--95, 1994.

\bibitem{cher95b}
I.~Cherednik.
\newblock Non-symmetric {Macdonald} polynomials.
\newblock {\em Int. Math. Res. Not.}, 10:483--515, 1995.

\bibitem{dunkl89a}
C.~F. Dunkl.
\newblock Differential-difference operators associated to reflection groups.
\newblock {\em Trans. Amer. Math. Soc.}, 311:167--183, 1989.


\bibitem{kakei96}
S.~Kakei.
\newblock Common algebraic structure for the {Calogero-Sutherland} models.
\newblock {\em J. Phys. A}, 29:L619--624, 1996.

\bibitem{kakei96b}
S.~Kakei.
\newblock An orthogonal basis for the {$B_N$}-type {Calogero} model.
\newblock {\em J. Phys. A}, 30:L535--541, 1997.

\bibitem{kakei97a}
S.~Kakei.
\newblock Intertwining operators for a degenerate double affine {Hecke}
algebra
and multivariable orthogonal polynomials.
\newblock q-alg/9706019.

\bibitem{kaneko93a}
J.~Kaneko.
\newblock Selberg integrals and hypergeometric functions associated with {Jack}
  polynomials.
\newblock {\em SIAM J. Math. Anal.}, 24:1086--1110, 1993.

\bibitem{kaneko96a}
J.~Kaneko.
\newblock {$q$-Selberg} integrals and {Macdonald} polynomials.
\newblock {\em Ann. Sci. \'Ec. Norm. Sup. $4^e$ s\'erie}, 29:1086--1110, 1996.

\bibitem{knop96a}
F.~Knop.
\newblock Integrality of two variable {Kostka} functions.
\newblock {\em J. Reine Ang. Math.}, 482:177--189, 1997.

\bibitem{knop96d}
F.~Knop and S.~Sahi.
\newblock Difference equations and symmetric polynomials defined by their
  zeros.
\newblock {\em Int. Math. Res. Not.}, 10:473--486, 1996.

\bibitem{knop96c}
F.~Knop and S.~Sahi.
\newblock A recursion and combinatorial formula for {Jack} polynomials.
\newblock {\em Inv. Math.}, 128:9--22, 1997.

\bibitem{lass97a}
M.~Lassalle.
\newblock Coefficients binomiaux g\'en\'eralis\'es et polyn\^omes de
  {Macdonald}.
\newblock preprint.

\bibitem{macunp1}
I.~G. Macdonald.
\newblock Hypergeometric functions.
\newblock Unpublished manuscript.

\bibitem{mac95}
I.~G. Macdonald.
\newblock Affine {Hecke} algebras and orthogonal polynomials.
\newblock {\em S\'eminaire Bourbaki, 47\`eme ann\'ee, Publ. I. R. M. A.
  Strasbourg}, 797, 1994-95.

\bibitem{mac}
I.~G. Macdonald.
\newblock {\em Symmetric functions and {Hall} polynomials}.
\newblock Oxford University Press, Oxford, 2nd edition, 1995.

\bibitem{noumi96c}
K.~Mimachi and M.~Noumi.
\newblock A reproducing kernel for nonsymmetric {Macdonald} polynomials.
\newblock q-alg/9610014.

\bibitem{okoun97a}
A.~Okounkov.
\newblock Binomial formula for {Macdonald} polynomials.
\newblock q-alg/9608021, to appear in {\it Math. Res. Lett.}

\bibitem{okoun96b}
A.~Okounkov.
\newblock Shifted {Macdonald} polynomials: q-integral representation and
  combinatorial formula.
\newblock q-alg/9605013, to appear in {\it Comp. Math.}

\bibitem{opdam95a}
E.~M. Opdam.
\newblock Harmonic analysis for certain representations of graded {Hecke}
  algebras.
\newblock {\em Acta Math.}, 175:75--121, 1995.

\bibitem{roesler97a}
M.~R\"osler.
\newblock Generalized {Hermite} polynomials and the heat equation for {Dunkl}
  operators.
\newblock q-alg/9703006, to appear in {\it Comm. Math. Phys.}

\bibitem{sahi97a}
S.~Sahi.
\newblock The binomial formula for nonsymmetric {Macdonald} polynomials.
\newblock q-alg/9703024.

\bibitem{sahi96a}
S.~Sahi.
\newblock A new scalar product for nonsymmetric {Jack} polynomials.
\newblock {\em Int. Math. Res. Not.}, 20:997--1004, 1996.

\bibitem{takemura97a}
K.~Takemura.
\newblock The {Yangian} symmetry in the spin {Calogero} model and its
  applications.
\newblock {\em J. Phys. A}, 30:6185--6204, 1997.

\bibitem{uji96a}
H.~Ujino and M.~Wadati.
\newblock Algebraic construction of the eigenstates for the second conserved
  operator of the quantum {Calogero} model.
\newblock {\em J. Phys. Soc. Japan}, 65:653--656, 1996.

\bibitem{uji96b}
H.~Ujino and M.~Wadati.
\newblock Rodrigues formula for Hi-Jack polynomials associated with
the quantum Calogero model.
\newblock {\em J. Phys. Soc. Japan}, 65:2423--2439, 1996.

\bibitem{uji97a}
H.~Ujino and M.~Wadati.
\newblock Orthogonality of the Hi-Jack polynomials associated with the
Calogero model.
\newblock {\em J. Phys. Soc. Japan}, 66:345--350, 1997.


\end{thebibliography}

\end{document}